*Research Article*

# The Solution to the BCS Gap Equation for Superconductivity and Its Temperature Dependence

**Shuji Watanabe**

*Division of Mathematical Sciences, Graduate School of Engineering, Gunma University, 4-2 Aramaki-machi, Maebashi 371-8510, Japan*

Correspondence should be addressed to Shuji Watanabe; shuwatanabe@gunma-u.ac.jp

From the viewpoint of operator theory, we deal with the temperature dependence of the solution to the BCS gap equation for superconductivity. When the potential is a positive constant, the BCS gap equation reduces to the simple gap equation. We first show that there is a unique nonnegative solution to the simple gap equation, that it is continuous and strictly decreasing, and that it is of class $C^2$ with respect to the temperature. We next deal with the case where the potential is not a constant but a function. When the potential is not a constant, we give another proof of the existence and uniqueness of the solution to the BCS gap equation, and show how the solution varies with the temperature. We finally show that the solution to the BCS gap equation is indeed continuous with respect to both the temperature and the energy under a certain condition when the potential is not a constant.

## 1. Introduction

We use the unit $k_B = 1$, where $k_B$ stands for the Boltzmann constant. Let $\omega_D > 0$ and $k \in \mathbb{R}^3$ stand for the Debye frequency and the wave vector of an electron, respectively. Let $h > 0$ be Planck's constant, and set $\hbar = h/(2\pi)$. Let $m > 0$ and $\mu > 0$ stand for the electron mass and the chemical potential, respectively. We denote by $T (\geq 0)$ the absolute temperature, and by $x$ the kinetic energy of an electron minus the chemical potential; that is, $x = \hbar^2 |k|^2/(2m) - \mu$. Note that $0 < \hbar\omega_D \ll \mu$.

In the BCS model [1, 2] of superconductivity, the solution to the BCS gap equation (1) is called the gap function. The gap function corresponds to the energy gap between the superconducting ground state and the superconducting first excited state. Accordingly, the value of the gap function (the solution) is nonnegative. We regard the gap function as a function of both $T$ and $x$ and denote it by $u$; that is, $u : (T, x) \mapsto u(T, x) (\geq 0)$. The BCS gap equation is the following nonlinear integral equation ($0 < \varepsilon \leq x \leq \hbar\omega_D$):

$$u(T, x) = \int_\varepsilon^{\hbar\omega_D} \frac{U(x, \xi) u(T, \xi)}{\sqrt{\xi^2 + u(T, \xi)^2}} \\ \times \tanh \frac{\sqrt{\xi^2 + u(T, \xi)^2}}{2T} d\xi, \quad (1)$$

where $U(\cdot, \cdot) > 0$ is the potential multiplied by the density of states per unit energy at the Fermi surface and is a function of $x$ and $\xi$. In (1) we introduce $\varepsilon > 0$, which is small enough and fixed ($0 < \varepsilon \ll \hbar\omega_D$). In the original BCS model, the integration interval is $[0, \hbar\omega_D]$; it is not $[\varepsilon, \hbar\omega_D]$. However, we introduce very small $\varepsilon > 0$ for the following mathematical reasons. In order to show the continuity of the solution to the BCS gap equation with respect to the temperature and in order to show that the transition to a superconducting state is a second-order phase transition, we make the form of the BCS gap equation somewhat easier to handle. So we choose the closed interval $[\varepsilon, \hbar\omega_D]$ as the integration interval in (1).

The integral with respect to $\xi$ in (1) is sometimes replaced by the integral over $\mathbb{R}^3$ with respect to the wave vector $k$. Odeh [3] and Billard and Fano [4] established the existence and uniqueness of the positive solution to the BCS gap equation in the case $T = 0$. For $T \geq 0$, Vansevenant [5] determined the transition temperature (the critical temperature) and showed that there is a unique positive solution to the BCS gap equation. Recently, Frank et al. [6] gave a rigorous analysis of the asymptotic behavior of the transition temperature at weak coupling. Hainzl et al. [7] proved that the existence of a positive solution to the BCS gap equation is equivalent to the existence of a negative eigenvalue of a certain linear operator to show the existence of a transition temperature. Moreover,

Hainzl and Seiringer [8] derived upper and lower bounds on the transition temperature and the energy gap for the BCS gap equation.

Since the existence and uniqueness of the solution were established for each fixed $T$ in the previous literature, the temperature dependence of the solution is not covered. It is well known that studying the temperature dependence of the solution to the BCS gap equation is very important in condensed matter physics. This is because, by dealing with the thermodynamical potential, this study leads to a mathematical proof of the statement that the transition to a superconducting state is a second-order phase transition. So, in condensed matter physics, it is highly desirable to study the temperature dependence of the solution to the BCS gap equation.

When the potential $U(\cdot,\cdot)$ in (1) is a positive constant, the BCS gap equation reduces to the simple gap equation (3). In this case, one assumes in the BCS model that there is a unique nonnegative solution to the simple gap equation (3) and that the solution is of class $C^2$ with respect to the temperature $T$ (see e.g., [1] and [9, (11.45), page 392]). In this paper, applying the implicit function theorem, we first show that this assumption of the BCS model indeed holds true; we show that there is a unique nonnegative solution to the simple gap equation (3) and that the solution is of class $C^2$ with respect to the temperature $T$. We next deal with the case where the potential is not a constant but a function. In order to show how the solution varies with the temperature, we then give another proof of the existence and uniqueness of the solution to the BCS gap equation (1) when the potential is not a constant. More precisely, we show that the solution belongs to the subset $V_T$ (see (12)). Note that the subset $V_T$ depends on $T$. We finally show that the solution to the BCS gap equation (1) is indeed continuous with respect to both $T$ and $x$ when $T$ satisfies (20) when the potential is not a constant.

Let

$$U(x,\xi) = U_1 \quad \text{at all } (x,\xi) \in [\varepsilon, \hbar\omega_D]^2, \quad (2)$$

where $U_1 > 0$ is a constant. Then the gap function depends on the temperature $T$ only. So we denote the gap function by $\Delta_1$ in this case; that is, $\Delta_1 : T \mapsto \Delta_1(T)$. Then (1) leads to the simple gap equation

$$1 = U_1 \int_\varepsilon^{\hbar\omega_D} \frac{1}{\sqrt{\xi^2 + \Delta_1(T)^2}} \tanh \frac{\sqrt{\xi^2 + \Delta_1(T)^2}}{2T} d\xi. \quad (3)$$

The following is the definition of the temperature $\tau_1 > 0$.

*Definition 1* (see [1]). Consider

$$1 = U_1 \int_\varepsilon^{\hbar\omega_D} \frac{1}{\xi} \tanh \frac{\xi}{2\tau_1} d\xi. \quad (4)$$

## 2. The Simple Gap Equation

Set

$$\Delta = \frac{\sqrt{(\hbar\omega_D - \varepsilon e^{1/U_1})(\hbar\omega_D - \varepsilon e^{-1/U_1})}}{\sinh(1/U_1)}. \quad (5)$$

**Proposition 2** (see [10, Proposition 2.2]). *Let $\Delta$ be as in* (5). *Then there is a unique nonnegative solution $\Delta_1 : [0, \tau_1] \to [0, \infty)$ to the simple gap equation* (3) *such that the solution $\Delta_1$ is continuous and strictly decreasing on the closed interval $[0, \tau_1]$:*

$$\Delta_1(0) = \Delta > \Delta_1(T_1) > \Delta_1(T_2)$$
$$> \Delta_1(\tau_1) = 0, \quad 0 < T_1 < T_2 < \tau_1. \quad (6)$$

*Moreover, the solution $\Delta_1$ is of class $C^2$ on the interval $[0, \tau_1)$ and satisfies*

$$\Delta_1'(0) = \Delta_1''(0) = 0, \quad \lim_{T\uparrow\tau_1}\Delta_1'(T) = -\infty. \quad (7)$$

*Proof.* Setting $Y = \Delta_1(T)^2$ turns (3) into

$$1 = U_1 \int_\varepsilon^{\hbar\omega_D} \frac{1}{\sqrt{\xi^2 + Y}} \tanh \frac{\sqrt{\xi^2 + Y}}{2T} d\xi. \quad (8)$$

Note that the right side is a function of the two variables $T$ and $Y$. We see that there is a unique function $T \mapsto Y$ defined by (8) implicitly, that the function $T \mapsto Y$ is continuous and strictly decreasing on $[0, \tau_1]$, and that $Y = 0$ at $T = \tau_1$. We moreover see that the function $T \mapsto Y$ is of class $C^2$ on the closed interval $[0, \tau_1]$. □

*Remark 3.* We set $\Delta_1(T) = 0$ for $T > \tau_1$.

*Remark 4.* In Proposition 2, $\Delta_1(T)$ is nothing but $\sqrt{f(T)}$ in [10, Proposition 2.2].

We introduce another positive constant $U_2 > 0$. Let $0 < U_1 < U_2$. We assume the following condition on $U(\cdot,\cdot)$:

$$U_1 \leq U(x,\xi)$$
$$\leq U_2 \quad \text{at all } (x,\xi) \in [\varepsilon, \hbar\omega_D]^2, \ U(\cdot,\cdot) \in C\left([\varepsilon, \hbar\omega_D]^2\right). \quad (9)$$

When $U(x,\xi) = U_2$ at all $(x,\xi) \in [\varepsilon, \hbar\omega_D]^2$, an argument similar to that in Proposition 2 gives that there is a unique nonnegative solution $\Delta_2 : [0, \tau_2] \to [0, \infty)$ to the simple gap equation

$$1 = U_2 \int_\varepsilon^{\hbar\omega_D} \frac{1}{\sqrt{\xi^2 + \Delta_2(T)^2}}$$
$$\times \tanh \frac{\sqrt{\xi^2 + \Delta_2(T)^2}}{2T} d\xi, \quad 0 \leq T \leq \tau_2. \quad (10)$$

Here, $\tau_2 > 0$ is defined by

$$1 = U_2 \int_\varepsilon^{\hbar\omega_D} \frac{1}{\xi} \tanh \frac{\xi}{2\tau_2} d\xi. \quad (11)$$

We again set $\Delta_2(T) = 0$ for $T > \tau_2$. A straightforward calculation gives the following.

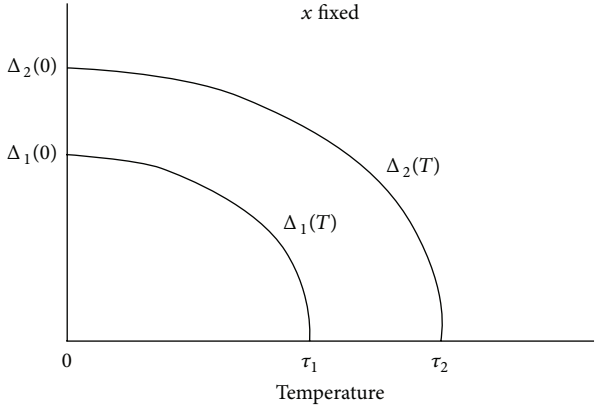

FIGURE 1: The graphs of the functions $\Delta_1$ and $\Delta_2$.

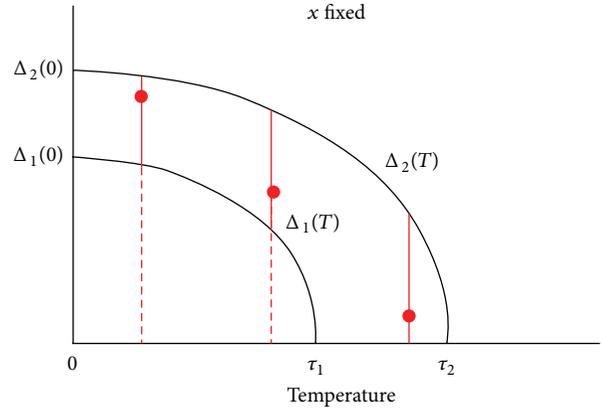

FIGURE 2: For each $T$, the solution $u_0(T, x)$ lies between $\Delta_1(T)$ and $\Delta_2(T)$.

**Lemma 5** ([11, Lemma 1.5]). *(a) The inequality $\tau_1 < \tau_2$ holds. (b) If $0 \leq T < \tau_2$, then $\Delta_1(T) < \Delta_2(T)$. If $T \geq \tau_2$, then $\Delta_1(T) = \Delta_2(T) = 0$.*

Note that Proposition 2 and Lemma 5 point out how $\Delta_1$ and $\Delta_2$ depend on the temperature and how $\Delta_1$ and $\Delta_2$ vary with the temperature; see Figure 1.

*Remark 6.* On the basis of Proposition 2, the present author [10, Theorem 2.3] proved that the transition to a superconducting state is a second-order phase transition under the restriction (2).

## 3. The BCS Gap Equation

Let $0 \leq T \leq \tau_2$ and fix $T$, where $\tau_2$ is that in (11). We consider the Banach space $C[\varepsilon, \hbar\omega_D]$ consisting of continuous functions of $x$ only and deal with the following subset $V_T$:

$$V_T = \{u(T, \cdot) \in C[\varepsilon, \hbar\omega_D] : \Delta_1(T) \leq u(T, x) \leq \Delta_2(T) \text{ at } x \in [\varepsilon, \hbar\omega_D]\}. \quad (12)$$

*Remark 7.* The subset $V_T$ depends on $T$. So we denote each element of $V_T$ by $u(T, \cdot)$; see Figure 1.

As it is mentioned in the introduction, the existence and uniqueness of the solution to the BCS gap equation were established for each fixed $T$ in the previous literature, and the temperature dependence of the solution is not covered. We therefore give another proof of the existence and uniqueness of the solution to the BCS gap equation (1) so as to show how the solution varies with the temperature. More precisely, we show that the solution belongs to $V_T$.

**Theorem 8** (see [11, Theorem 2.2]). *Assume condition (9) on $U(\cdot, \cdot)$. Let $T \in [0, \tau_2]$ be fixed. Then there is a unique nonnegative solution $u_0(T, \cdot) \in V_T$ to the BCS gap equation (1) ($x \in [\varepsilon, \hbar\omega_D]$):*

$$u_0(T, x) = \int_\varepsilon^{\hbar\omega_D} \frac{U(x, \xi) u_0(T, \xi)}{\sqrt{\xi^2 + u_0(T, \xi)^2}} \times \tanh \frac{\sqrt{\xi^2 + u_0(T, \xi)^2}}{2T} d\xi. \quad (13)$$

*Consequently, the solution is continuous with respect to $x$ and varies with the temperature as follows:*

$$\Delta_1(T) \leq u_0(T, x) \leq \Delta_2(T) \quad at \ (T, x) \in [0, \tau_2] \times [\varepsilon, \hbar\omega_D]. \quad (14)$$

*Proof.* We define a nonlinear integral operator $A$ on $V_T$ by

$$Au(T, x) = \int_\varepsilon^{\hbar\omega_D} \frac{U(x, \xi) u(T, \xi)}{\sqrt{\xi^2 + u(T, \xi)^2}} \times \tanh \frac{\sqrt{\xi^2 + u(T, \xi)^2}}{2T} d\xi, \quad (15)$$

where $u(T, \cdot) \in V_T$. Clearly, $V_T$ is a bounded, closed, and convex subset of the Banach space $C[\varepsilon, \hbar\omega_D]$. A straightforward calculation gives that the operator $A : V_T \to V_T$ is compact. Therefore, the Schauder fixed point theorem applies, and hence the operator $A : V_T \to V_T$ has at least one fixed point $u_0(T, \cdot) \in V_T$. Moreover, we can show the uniqueness of the fixed point; see Figure 2. □

The existence of the transition temperature $T_c$ is pointed out in the previous papers [5–8]. In our case, it is defined as follows.

*Definition 9.* Let $u_0(T, \cdot) \in V_T$ be as in Theorem 8. The transition temperature $T_c$ stemming from the BCS gap equation (1) is defined by

$$T_c = \inf \{T > 0 : u_0(T, x) = 0 \text{ at all } x \in [\varepsilon, \hbar\omega_D]\}. \quad (16)$$

*Remark 10.* Combining Definition 9 with Theorem 8 implies that $\tau_1 \leq T_c \leq \tau_2$. For $T > T_c$, we set $u_0(T,x) = 0$ at all $x \in [\varepsilon, \hbar\omega_D]$.

## 4. Continuity of the Solution with respect to the Temperature

Let $U_0 > 0$ be a constant satisfying $U_0 < U_1 < U_2$. An argument similar to that in Proposition 2 gives that there is a unique nonnegative solution $\Delta_0 : [0, \tau_0] \to [0, \infty)$ to the simple gap equation

$$1 = U_0 \int_\varepsilon^{\hbar\omega_D} \frac{1}{\sqrt{\xi^2 + \Delta_0(T)^2}}$$

$$\times \tanh \frac{\sqrt{\xi^2 + \Delta_0(T)^2}}{2T} d\xi, \quad 0 \leq T \leq \tau_0. \tag{17}$$

Here, $\tau_0 > 0$ is defined by

$$1 = U_0 \int_\varepsilon^{\hbar\omega_D} \frac{1}{\xi} \tanh \frac{\xi}{2\tau_0} d\xi. \tag{18}$$

We set $\Delta_0(T) = 0$ for $T > \tau_0$. A straightforward calculation gives the following.

**Lemma 11.** *(a)* $\tau_0 < \tau_1 < \tau_2$.
*(b)* If $0 \leq T < \tau_0$, then $0 < \Delta_0(T) < \Delta_1(T) < \Delta_2(T)$.
*(c)* If $\tau_0 \leq T < \tau_1$, then $0 = \Delta_0(T) < \Delta_1(T) < \Delta_2(T)$.
*(d)* If $\tau_1 \leq T < \tau_2$, then $0 = \Delta_0(T) = \Delta_1(T) < \Delta_2(T)$.
*(e)* If $\tau_2 \leq T$, then $0 = \Delta_0(T) = \Delta_1(T) = \Delta_2(T)$.

*Remark 12.* Let the functions $\Delta_k$ ($k = 0, 1, 2$) be as above. For each $\Delta_k$, there is the inverse $\Delta_k^{-1} : [0, \Delta_k(0)] \to [0, \tau_k]$. Here,

$$\Delta_k(0) = \frac{\sqrt{(\hbar\omega_D - \varepsilon e^{1/U_k})(\hbar\omega_D - \varepsilon e^{-1/U_k})}}{\sinh(1/U_k)}, \tag{19}$$

and $\Delta_0(0) < \Delta_1(0) < \Delta_2(0)$.

We introduce another temperature. Let $T_1$ satisfy $0 < T_1 < \Delta_0^{-1}(\Delta_0(0)/2)$ and

$$\frac{\Delta_0(0)}{4\Delta_2^{-1}(\Delta_0(T_1))} \tanh \frac{\Delta_0(0)}{4\Delta_2^{-1}(\Delta_0(T_1))}$$

$$> \frac{1}{2}\left(1 + \frac{4\hbar^2\omega_D^2}{\Delta_0(0)^2}\right). \tag{20}$$

*Remark 13.* Numerically, the temperature $T_1$ is very small.

Consider the following subset $V$ of the Banach space $C([0, T_1] \times [\varepsilon, \hbar\omega_D])$ consisting of continuous functions of both the temperature $T$ and the energy $x$:

$$V = \{u \in C([0, T_1] \times [\varepsilon, \hbar\omega_D]) : \Delta_1(T) \leq u(T,x)$$

$$\leq \Delta_2(T) \text{ at } (T,x) \in [0, T_1] \times [\varepsilon, \hbar\omega_D]\}. \tag{21}$$

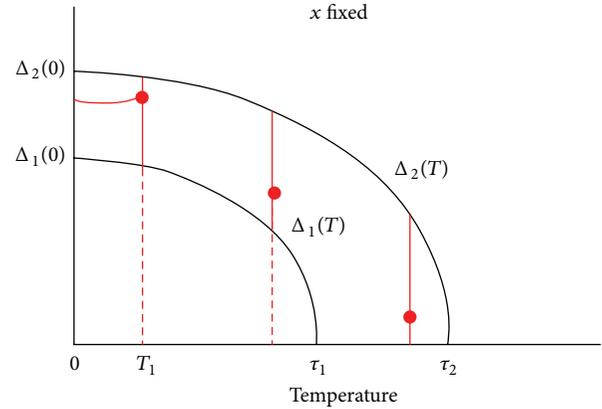

FIGURE 3: The solution $u_0$ is continuous on $[0, T_1] \times [\varepsilon, \hbar\omega_D]$.

**Theorem 14** (see [12, Theorem 1.2]). *Assume* (9). *Let $u_0$ be as in Theorem 8 and $V$ as in* (21). *Then $u_0 \in V$. Consequently, the gap function $u_0$ is continuous on $[0, T_1] \times [\varepsilon, \hbar\omega_D]$.*

*Proof.* We define a nonlinear integral operator $B$ on $V$ by

$$Bu(T,x) = \int_\varepsilon^{\hbar\omega_D} \frac{U(x,\xi) u(T,\xi)}{\sqrt{\xi^2 + u(T,\xi)^2}}$$

$$\times \tanh \frac{\sqrt{\xi^2 + u(T,\xi)^2}}{2T} d\xi, \tag{22}$$

where $u \in V$.

Clearly, $V$ is a closed subset of the Banach space $C([0, T_1] \times [\varepsilon, \hbar\omega_D])$. A straightforward calculation gives that the operator $B : V \to V$ is contractive as long as (20) holds true. Therefore, the Banach fixed-point theorem applies, and hence the operator $B : V \to V$ has a unique fixed point $u_0 \in V$. The solution $u_0 \in V$ to the BCS gap equation is thus continuous both with respect to the temperature and with respect to the energy $x$; see Figure 3. □

## Acknowledgment

Shuji Watanabe is supported in part by the JSPS Grant-in-Aid for Scientific Research (C) 24540112.